\documentclass[12pt]{article}
 \usepackage{amsfonts}
\usepackage{theorem}
\usepackage[OT2,OT1]{fontenc}
\newcommand\cyr{%
\renewcommand\rmdefault{wncyr}%
\renewcommand\sfdefault{wncyss}%
\renewcommand\encodingdefault{OT2}%
\normalfont
\selectfont}
\DeclareTextFontCommand{\textcyr}{\cyr}

\newtheorem{theorem}{Theorem}[section]

\newtheorem{lemma}{Lemma}[section]
\newtheorem{proposition}{Proposition}[section]

\newtheorem{example}{Example}
\newtheorem{constr}{Construction}[section]
\usepackage[russian, english]{babel}
\usepackage{centernot}
\usepackage{amsmath}
\usepackage[unicode,colorlinks,plainpages=false]{hyperref}
\newcommand{\openbox}{$\begin{array}{c}
\hspace*{-0.55em}\sqcap \hspace*{-0.60em}\\[-0.4em] \hline
\multicolumn{1}{c}{\hspace*{-0.60em}}\\[-0.8em]
\end{array}$}

\begin{document}

\centerline{\bf On the Probability}
\centerline{\bf That Two Elements of a Finite Semigroup}
\centerline{\bf Have the Same Right Matrix \footnote{This work was supported by the National Research, Development and Innovation Office – NKFIH, 115288}}

\bigskip

\centerline{\bf Attila Nagy and Csaba T\'oth}

\bigskip

\begin{abstract}In this paper we study the probability that two elements selected at random with replacement from a given finite semigroup act the same by right translation on the semigroup, that is, the chosen elements have the same right matrix.
\end{abstract}

\bigskip


\medskip


\section{Introduction and motivation}

There are many papers in the literature of mathematics that use probabilistic methods in the study of special algebraic structures. Here we only cite papers \cite{Ahmedidelir}, \cite{Dixon1}, \cite{Dixon2}, \cite{Eberhard}, \cite{Gustafson} and \cite{MacHale}, because we will only refer to them. All of these papers
%Papers \cite{Ahmedidelir}, \cite{Dixon1}, \cite{Dixon2}, \cite{Eberhard}, \cite{Gustafson} and \cite{MacHale}
deal with special cases of the following problem: for a given algebraic structure $A$ and a given binary relation $\sigma$ on $A$, find the probability $P_{\sigma}(A)$ that $(a,b)\in \sigma$ is satisfied for two elements $a$ and $b$ of $A$ selected at random with replacement.
We note that random elements are chosen independently with the uniform
distribution on $A$. Thus every pair $(a, b)\in A\times A$ has the same probability $\frac{1}{|A|^2}$ of being chosen and so
\[P_{\sigma}(A)=\frac{|\{(a,b)\in A\times A:\ (a,b)\in \sigma\}|}{|A|^2}.\]

In \cite{Ahmedidelir}, \cite{Gustafson} and \cite{MacHale} the probability $P_{\sigma}(A)$ is examined in that cases, when $A$ is a finite non-commutative semigroup, a non-commutative group and a non-commutative ring, respectively; in all three cases $\sigma$ is defined by $(x, y)\in \sigma$ for $x, y\in A$ if and only if $xy=yx$.
In \cite{Dixon2}, the probability $P_{\sigma}(A)$ is investigated in that case when $A$ is a finite simple group and $\sigma$ is defined by $(x, y)\in \sigma$ for $x, y\in A$ if and only if $w(x,y)=e$, where $w$ is a given non-trivial element of the free group $F_2$ and $e$ is the identity element of the group $A$. In both of \cite{Dixon1} and \cite{Eberhard}, the probability $P_{\sigma}(A)$ is examined in that case when $A$ is the symmetric group $S_n$ of degree $n$. In \cite{Dixon1}, $\sigma$ is defined by $(x,y)\in \sigma$ for $x, y\in A$ if and only if $x$ and $y$ generate $S_n$. In \cite{Eberhard}, $\sigma$ is defined by $(x,y)\in \sigma$ for $x, y\in A$ if and only if $x$ and $y$ generate the alternating group $A_n$ or the symmetric group $S_n$.

\medskip
The above mentioned investigations motivate us to examine the probability $P_{\sigma }(A)$ in a further special case. In our present paper we deal with the probability $P_{\theta _S}(S)$, where $S$ is a finite semigroup and $\theta _S$ is the kernel of the right regular matrix representation of $S$. With other words, $P_{\theta _S}(S)$ is the
probability that two elements $a$ and $b$ of a given finite semigroup $S$ selected at random with replacement act the same by right translation on the semigroup $S$, that is, $xa=xb$ is satisfied for all $x\in S$.

\section{Preliminaries}

Let $S$ be a semigroup and $G^0$ be a semigroup arising from a one-element group $G=\{ 1\}$ by adjunction of a zero element $0$. By an $S\times S$ matrix over $G^0$ we mean a mapping of $S\times S$ into $G^0$. Let $A$ be an $S\times S$ matrix over $G^0$. For an element $s\in S$, the set $\{ A(s, x):\ x\in S\}$ is called the $s-row$ of $A$. An $S\times S$ matrix $A$ over $G^0$ is called strictly row-monomial if each row of $A$ contains exactly one non-zero element of $G^0$.

For an element $a$ of a semigroup $S$, let $R^{(a)}$ denote the strictly row-monomial $S$-matrix over $G^0$ defined by
\[ {R}^{(a)}((x,y))=\begin{cases}
1, & \text{if $xa=y$}\\
0 & \text{otherwise.}
\end{cases}\]

This matrix is called the right matrix over $G^0$ defined by $a$. (The right matrices were investigated, for example, in \cite{Nagy0:sg-7} and \cite{NagyRonyai:sg-8} in that case when $1$ is the identity element and $0$ is the zero element of a field.)
It is known that \[a \mapsto R^{(a)}\] is a homomorphism of a semigroup $S$ into the multiplicative semigroup of all strictly row-monomial $S$-matrices over $G^0$ (see \cite[Exercise 4(b) for \S 3.5]{Clifford1:sg-1}).
This homomorphism is called the right regular matrix representation of the semigroup $S$. Let $\theta _S$ denote the kernel of the right regular matrix representation of $S$. It is clear that
\[\theta _S=\{ (a, b)\in S\times S:\ (\forall x\in S)\ xa=xb\}.\]

\medskip

In this paper, we investigate the probability that two elements of a finite semigroup $S$ selected at random with replacement have the same right matrix,
that is, we investigate the probability \[P_{\theta _S}(S)=\frac{|\{(a, b)\in S\times S: (a, b)\in \theta _S\}}{|S|^2}.\]
For brevity, $P_{\theta _S}(S)$ will be denoted by $P_{\theta}(S)$.

In our investigation, the following construction (which is a special case of the construction in \cite[Theorem 1]{NagyKolib}) will play an important role.

\begin{constr}\rm (\cite[Construction 1]{Acta})\label{constr1} Let $T$ be a left cancellative semigroup (that is, a semigroup with the property that $xa=xb$ implies $a=b$ for every $x, a, b\in T$). For each $t\in T$, associate a nonempty set $S_t$ such that $S_t\cap S_r=\emptyset$ for all $t, r\in T$ with $t\neq r$. As $T$ is left cancellative, $x \mapsto tx$ is an injective mapping of $T$ onto $tT$.

For arbitrary couple $(t, r)\in T\times T$ with $r\in tT$, let $(\cdot )\varphi _{t, r}$ be a mapping of $S_t$ into $S_r$ acting on the right. For all $t\in T$, $r\in tT$, $q\in rT\subseteq tT$
and $a\in S_t$, assume
\begin{equation}\label{mappings}
(a)(\varphi _{t,r}\circ \varphi _{r,q})=(a)\varphi _{t, q}.
\end{equation}

On the set $S=\cup _{t\in T}S_t$ define an operation $\star$ as follows: for arbitrary $a\in S_t$ and $b\in S_x$, let \[a\star b=(a)\varphi _{t, tx}.\]

By \cite[Theorem 1]{NagyKolib}, $(S; \star )$ is a semigroup, and each set $S_t$ ($t\in T$) is a $\theta _S$-class of $S$.\hfill\openbox
\end{constr}

\medskip

For notions not defined in this paper, we refer to the books \cite{Clifford1:sg-1} and \cite{Nagybook:sg-5}.

\section{The results}

Let $A$ be a non-empty set and $\sigma$ be a binary relation on $A$. Let $P_{\sigma}(A)$ denote the probability that $(a,b)\in \sigma$ is satisfied for two elements $a$ and $b$ of $A$ selected at random with replacement.

\begin{theorem}\label{thm1} Let $p$ be an arbitrary rational number with $0\leq p\leq 1$. Then the following assertions are equivalent.
\begin{itemize}
\item[(i)] There is a finite semigroup $S$ such that $P_{\theta }(S)=p$.
\item[(ii)] There is a non-empty finite set $A$ and an equivalence relation $\sigma$ on $A$ such that $P_{\sigma}(A)=p$.
\end{itemize}
\end{theorem}

\noindent
{\bf Proof}. We need only to prove that $(ii)$ implies $(i)$. Assume $(ii)$. Let $T$ be a commutative group of order $|A/\sigma|$.
%Define an operation on $T$: $tr=r$ for every $t, r\in T$. Then $T$ forms a right zero semigroup under this operation.
Let $S_t$ ($t\in T$) denote the $\sigma$-classes of $A$. For every $t\in T$, fix an element $s_t$ in $S_t$. For every $t, r\in T$, let $(\cdot )\varphi _{t, r}$ be the mapping of $S_t$ into $S_r$ which maps the elements of $S_t$ to $s_r$. It is easy to see that the family $\{\varphi _{t, r}:\ t, r\in T\}$ of mappings satisfies the following condition: for every $t, r\in T$ and every $a\in S_t$,
\[(a)(\varphi _{t,r}\circ \varphi _{r,q})=(a)\varphi _{t,q}.\] Thus condition (\ref{mappings}) of Construction~\ref{constr1} is satisfied. Hence $S=\cup _{t\in T}S_t$ forms a semigroup under the operation $\star$ defined by \[a\star b=(a)\varphi _{t, tr}=s_{tr}\] for every $a\in S_t$ and $b\in S_r$. As
\[a\star b=(a)\varphi _{t, tr}=s_{tr}=s_{rt}=(b)_{r, rt}=b\star a\] for every $r, t\in T$ and $a\in S_t$, $b\in S_r$, the semigroup $S$ is commutative. As $T$ is left reductive, \cite[Theorem 1]{NagyKolib} implies that the $\theta$-classes of the semigroup $S$ are the sets $S_t$ ($t\in T$). Consequently $P_{\theta}(S)=P_{\sigma}(A)=p$.\hfill\openbox

\bigskip

By Theorem~\ref{thm1}, we can concentrate our attention on non-empty finite sets $A$ and equivalence relations $\sigma$ on $A$.

\begin{theorem} \label{thm2} If $\sigma$ is an equivalence relation on a non-empty finite set $A$, then $P_{\sigma}(A)\geq \frac{1}{|A/\sigma|}$. The equation $P_{\sigma}(A)=\frac{1}{|A/\sigma|}$ holds if and only if the cardinality of any two $\sigma$-classes of $A$ is the same.
\end{theorem}

\noindent {\bf Proof}. Let $m$ denote the cardinality of the factor set $A/\sigma$. Let $A_i$ ($i=1,\dots ,m$) denote the pairwise different $\sigma$-classes of $A$. If $|A_i|=t_i$ then
\[P_{\sigma}(A)=\frac{t_1^2+\cdots +t_m^2}{(t_1+\cdots +t_m)^2}.\]
By the well known connection between the root mean square and the arithmetic mean, we have
\[\sqrt{\frac{\sum_{i=1}^{m}t_i^2}{m}}\geq \frac{\sum_{i=1}^{m}t_i}{m},\] that is
\[\frac{\sum_{i=1}^{m}t_i^2}{m}\geq \frac{(\sum_{i=1}^{m}t_i)^2}{m^2}\]
from which we get
\[P_{\sigma}(A)=\frac{t_1^2+\cdots +t_m^2}{(t_1+\cdots +t_m)^2}\geq \frac{1}{m}=\frac{1}{|A/\sigma|}.\]
The equation $P_{\sigma}(A)=\frac{1}{|A/\sigma|}$ holds if and only if
\[\frac{t_1^2+\cdots +t_m^2}{(t_1+\cdots +t_m)^2}=\frac{1}{m},\] that is
\[\sqrt{\frac{\sum_{i=1}^{m}t_i^2}{m}}=\frac{\sum_{i=1}^{m}t_i}{m}\] which is satisfied if and only if
\[t_1=\cdots =t_m.\]
Thus the theorem is proved.
\hfill\openbox

\medskip

A semigroup $S$ is said to be left reductive if, for every $a, b\in S$, the assumption "$xa=xb$ for all $x\in S$" implies $a=b$. It is clear that a semigroup $S$ is left reductive if and only if $\theta _S$ is the identity relation on $S$. The following proposition is a consequence of Theorem~\ref{thm1} and Theorem~\ref{thm2}.

\begin{proposition}\label{thm3} A finite semigroup $S$ is left reductive if and only if \[P_{\theta}(S)=\frac{1}{|S|}.\]\hfill\openbox
\end{proposition}

\medskip

The following two theorems show that, for non left reductive semigroups $S$, the probability $P_{\theta}(S)$ can be arbitrarily closed to $0$ and also to $1$.

\begin{theorem}\label{thm4} For any positive integer $n$, there is a non left reductive finite semigroup $S^{(n)}$ such that
$\lim _{n\to \infty} P_{\theta}(S^{(n)})=0$.
\end{theorem}

\noindent
{\bf Proof}. Let $A^{(n)}$ be a set and $\sigma$ is an equivalence relation on $A^{(n)}$ such that \[|A^{(n)}/\sigma|=n,\] and every $\sigma$-classes of $A$ contains $k\geq 2$ elements. By Theorem~\ref{thm2}, \[P_{\sigma}(A^{(n)})=\frac{1}{n}.\] Applying Theorem~\ref{thm1}, for every positive integer $n$, there is a semigroup $S^{(n)}$ such that \[P_{\theta}(S^{(n)})=P_{\sigma}(A^{(n)}).\] Hence \[\lim _{n\to \infty} P_{\theta}(S^{(n)})=0.\]\hfill\openbox

\begin{theorem}\label{thm73} For every integers $n$, there is a non left reductive finite semigroup $S^{(n)}$ such that
\[lim _{n\to \infty}P_{\theta}(S^{(n)})=1 .\]
\end{theorem}

\noindent
{\bf Proof}. Let $k_0\geq 2$ be a fixed integer. For arbitrary positive integer $n$, let $A^{(n)}$ be a set containing $n+k_0$ elements. Let $\sigma$ be an equivalence relation on $A$ with the following conditions: \[|A/\sigma|=k_0,\] moreover one of the $\sigma$-classes on $A$ contains $n+1$ elements and the other $\sigma$-classes on $A$ are singletons (that is, they contain one element). It is clear that
\[P_{\sigma}(A^{(n)})=\frac{(n+1)^2+k_0-1}{(n+k_0)^2}.\] It is a matter of checking to see that
\[lim_{n\to \infty}\frac{(n+1)^2+k_0-1}{(n+k_0)^2}=1.\] Hence
\[lim _{n\to \infty}P_{\theta}(A^{(n)})=1.\]
By Theorem~\ref{thm1}, for every positive integer $n$, there is a semigroup $S^{(n)}$ such that
\[P_{\sigma}(S^{(n)})=P_{\sigma}(A^{(n)}).\]
This proves our theorem.\hfill\openbox

\medskip

For an arbitrary finite semigroup $S$, $P_{\theta}(S)=1$ is satisfied if and only if $\theta _S=\omega _S$ ($\omega _S$ is the universal relation on $S$), that is, $xa=xb$ holds for every $x, a, b\in S$. The next theorem characterizes not necessarily finite semigroups $S$ in which $\theta _S$ is the universal relation on $S$.

We shall use the following notions. A homomorphism $\varphi$ of a semigroup $S$ onto an ideal $I$ of $S$ is called a retract homomorhism if $\varphi$ leaves the elements of $I$ fixed. An ideal $I$ of a semigroup $S$ is called a retract ideal if there is a retract homomorphism of $S$ onto $I$. In this case we say that $S$ is a retract (ideal) extension of $I$ by the Rees factor semigroup $S/I$. A semigroup satisfying the identity $ab=a$ is called a left zero semigroup. A semigroup with a zero $0$ is called a zero semigroup if it satisfies the identity $ab=0$.

\begin{theorem}\label{thmuniv}
For a semigroup $S$, $\theta _S$ is the universal relation on $S$ if and only if $S$ is a retract extension of a left zero semigroup by a zero semigroup.
\end{theorem}

\noindent
{\bf Proof}. Let $S$ be a semigroup in which $\theta _S$ is the universal relation on $S$. Then $xa=xb$ holds for every $x, a, b\in S$. Let $a\in S$ be an arbitrary element. Then \[a^2=aa=aa^2=a^3\] and so \[(a^2)^2=aa^3=aa^2=a^3=a^2,\] that is, $a^2$ is an idempotent element. Let $E(S)$ denote the set of all idempotent elements of $S$. As $ab=a^2\in E(S)$ for every $a, b\in S$, the set $E(S)$ is an ideal of $S$, and the Rees factor semigroup $Q=S/E(S)$ is a zero semigroup. For arbitrary $e, f\in E(S)$, \[ef=ee=e\] and so $E(S)$ is a left zero semigroup. For every $a\in S$, we have $aS\subseteq E(S)$ and $|aS|=1$. For every $a\in S$, let $\varphi (a)$ denote the element of $aS$. By the above, \[\varphi (a)\in E(S)\] for every $a\in S$. Moreover, $\varphi (e)=e$ for every idempotent element $e$ of $S$. Let $x^*\in S$ be an arbitrary fixed element. Then, for every $a, b\in S$, we have \[\varphi (ab)=abx^*=ax^*bx^*=\varphi (a)\varphi (b).\] Hence $\varphi$ is a homomorphism of $S$ onto $E(S)$. As $\varphi$ leaves the elements of $E(S)$ fixed, it is a retract homomorphism of $S$ onto $E(S)$. Thus $S$ is a retract extension of the left zero semigroup $E(S)$ by the zero semigroup $Q=S/E(S)$.

Conversely, let $S$ be a semigroup and $I$ is an ideal of $S$ such that $I$ is a left zero semigroup, the Rees factor semigroup $S/I$ is a zero semigroup, and there is a retract homomorphism $\varphi$ of $S$ onto $I$. Then, for arbitrary $x, a, b\in S$, we have $xa, xb\in I$ and so
\[xa=\varphi (xa)=\varphi (x)\varphi (a)=\varphi (x)=\varphi (x)\varphi (b)=\varphi (xb)=xb.\]
Hence $\theta _S$ is the universal relation on $S$. Thus the theorem is proved.
\hfill\openbox

\medskip

In the next we show how to construct semigroups $S$ in which $\theta _S$ is the universal relation on $S$.

\begin{constr}\label{construction2} Let $S$ be a non-empty set and $L$ is a non-empty subset of $S$. Let $(\cdot )\varphi$ be an arbitrary mapping of $S$ onto $L$ which leaves the elements of $L$ fixed. Define an operation $\star$ on $S$ as follows: for arbitrary $a, b\in S$, let $a\star b=(a)\varphi$.
For every $a, b, c\in S$, \[a\star (b\star c)=a\star (b)\varphi=(a)\varphi=((a)\varphi )\varphi=(a\star b)\varphi=(a\star b)\star c,\]
that is, $S$ is a semigroup with the operation $\star$. This semigroup is denoted by $(S, L, \varphi, \star )$.
\end{constr}

\begin{theorem}\label{thmconstr} In the semigroup $S=(S, L, \varphi, \star )$, the equation $\theta _S=\omega _S$ is satisfied. Conversely, every semigroup $S$ in which $\theta _S=\omega _S$ is satisfied is isomorphic to a semigroup defined in Construction~\ref{construction2}.
\end{theorem}

\noindent
{\bf Proof}.  For every $a, b\in (S, L, \varphi, \star )$, we have  $a\star b\in L$. Thus $L$ is an ideal of $S$ and the Rees factor semigroup $S/L$ is a zero semigroup. For every $a, b\in L$, we have
\[a\star b =(a)\varphi =a.\] Thus $L$ is a left zero semigroup. As \[(a\star b)\varphi =((a)\varphi)\varphi=(a)\varphi=((a)\varphi)((b)\varphi),\] $\varphi$ is a ratract homomorphism of $S$ onto $L$. Thus $S=(S, L, \varphi, \star )$ is a retract extension of the left zero semigroup $L$ by the zero semigroup $S/L$. Consequently $\theta _S=\omega _S$ by Theorem~\ref{thmuniv}.

Conversely, assume that $S$ is a semigroup in which $\theta _S=\omega _S$. By Theorem~\ref{thmuniv}, there is an ideal $L$ of $S$ such that $L$ is a left zero semigroup, the Rees factor semigroup $S/L$ is a zero semigroup, and there is a retract homomorphism $\varphi$ of $S$ onto $L$. Consider the semigroup $(S, L, \varphi, \star )$ defined as in Construction~\ref{construction2}. As
\[ab=(ab)\varphi =(a)\varphi (b)\varphi =(a)\varphi =a\star b\]
for every $a, b\in S$,
the semigroups $S$ and $(S, L, \varphi, \star )$ are isomorphic. \hfill\openbox

\bigskip

\noindent
Attila Nagy and Csaba T\'oth

\noindent
Department of Algebra, Institute of Mathematics,

\noindent
Budapest University of Technology and Economics,

\noindent
1111 Budapest, Egry J\'ozsef u. 1, Hungary

\noindent
e-mails: nagyat@math.bme.hu (A. Nagy); tcsaba94@gmail.com (Cs. T\'oth)

\end{document}